# Gitterwege reeller Zahlen

Lorenz Friess

Anstract
Let r be a real number, 0<r<1, given as a dual number. With the sequence of "0" "1" we define a path in a hexagonal lattice. Generalisations to other bases and lattices. Relations between the properties of r and its path are considered. Related questions are discussed.

Zusammenfassung
Sei 0<r<1 eine Dualzahl. Mit ihren Ziffern wird ein Weg in einem hexagonalen Gitter definiert. Welche Bezeihungen bestehen zwischen Eigenschaften von r und Eigenschaften des Weges?

Gegeben sei eine reelle Zahl r, 0 <= r <= 1. r sei zur Basis zwei dargestellt. Ihre Ziffern bilden eine Folge von "0" und "1". Diese Folge bezeichnen wir mit $(z_i)$.
Mit dieser Folge definieren wir einen Weg W = W(r) in einem hexagonalen Gitter beginnend bei einem (beliebig gewählten, festen) Startpunkt S mit der Anfangsrichtung "nach oben" (Abb.). Für eine "0" gehen wir nach links, für eine "1" nach rechts zum nächsten Punkt im Gitter und verfahren dort unter Verwendung der zweiten Ziffer ebenso usw.

Der Gitterpunkt, der mit $z_i$ erreicht wird, werde mit $G_i$ bezeichnet.

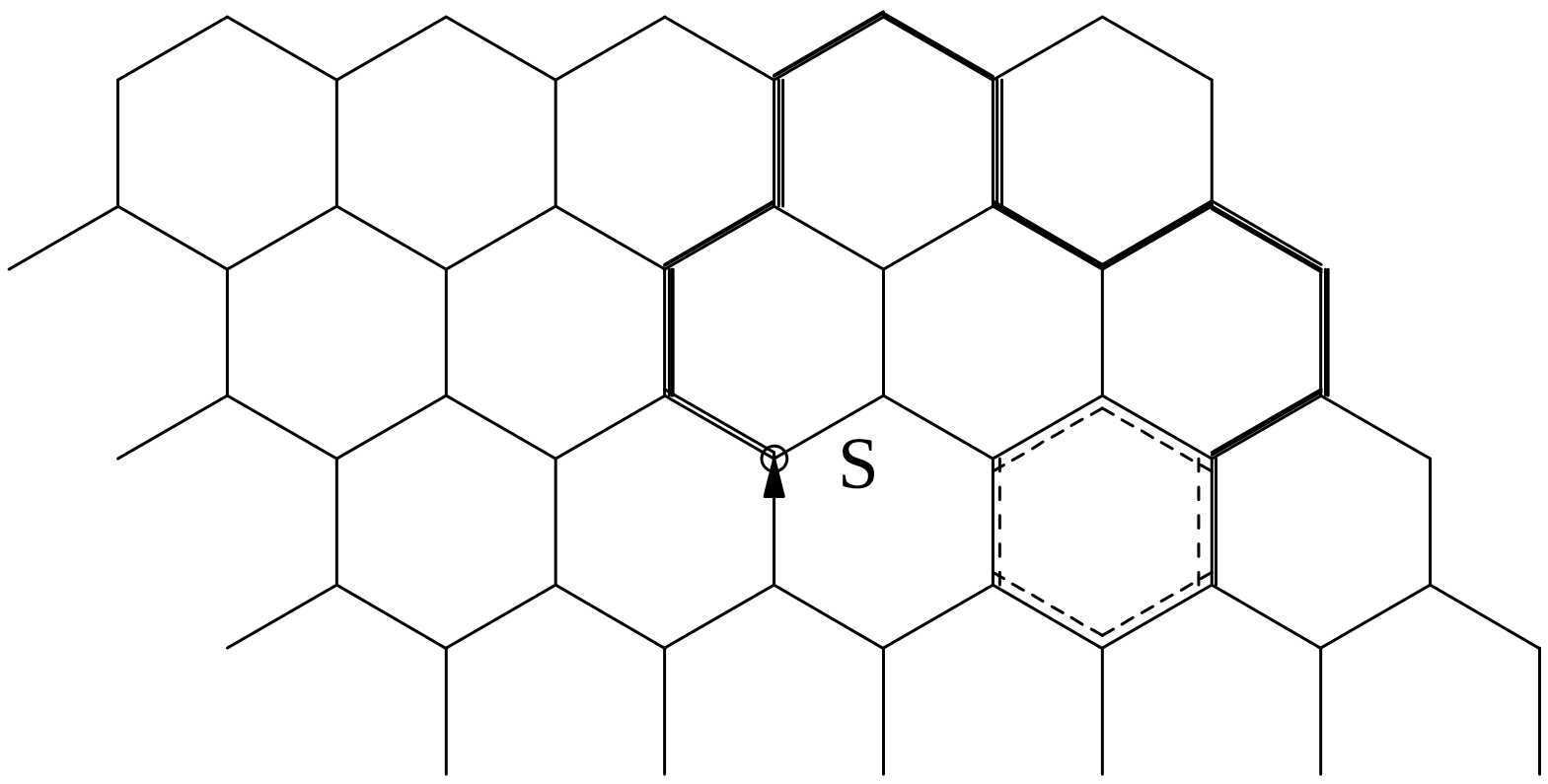

Für r = 0.0110111001111 ist W(r) in der Abbildung durch Strichverstärkung angedeutet. Sind alle folgenden Ziffern = "0", so liefern je sechs einen Durchlauf desselben Sechsecks (gestrichelt).

Wir interessieren uns u. a. für den Zusammenhang zwischen Eigenschaften des Weges einer Zahl r und Eigenschaften von r, z.B. rational, algebraisch oder transzendent zu sein.

Definition 1
Der Abstand zweier Punkte A, B im Gitter d(A, B) sei das Minimum der Länge aller Wege von A nach B. Die Länge eines Weges ist selbsterklärend. Der Abstand d(S, A) wird mit d(A) bezeichnet.

Es ist d(A, B) >= 0, d(A, B) = 0 <=> A = B, d(A, C) <= d(A, B) + d(B, C), d. h. d ist eine Metrik.

Definition 2

Ein Weg heißt beschränkt, wenn es ein M gibt mit d(A) < M für alle Punkte A des Weges.

Ist r rational besitzt ($z_i$) eine Periode. Der Teil des Weges, der den Ziffern der Periode entspricht, wird fortlaufend aneinandergesetzt. Sind die Richtungen am Anfang und Ende des Teilwegs gleich ist der Weg unbeschränkt. Sind die Richtungen am Anfang und Ende des Teilwegs verschieden bilden zwei, drei oder sechs Teilwege einen geschlossenen Weg.
Der Weg von r = 0.1010... ( = 2/3 ) ist unbeschränkt, der Weg von r = 0.11001100.. ( = 4/5 ) ebenso.
Der Weg von r = 0.110110... ( = 6/7) durchläuft immer wieder dieselben Gitterpunkte (es sind 18).

Definition 3 (Windungszahl)
Die Windungszahl w(W, A) sei die Anzahl wie oft der Weg W vom Anfang S bis zum Punkt A den Punkt S umläuft. (Wenn W durch S geht werde dies nicht mitgezählt.)
Allgemeiner sei w(W, A, B) die Anzahl wie oft W vom Anfang S bis zum Punkt A den Punkt B umläuft.

Definition 4 (Torsionszahl)
An jedem Gitterpunkt erfolgt eine Richtungsänderung um $\pi/3$. Es sei R(W, A) = Summe Richtungsänderungen von W von S bis A. Die Torsionszahl t(W, A) sei R(W, A)/6.

Wir betrachten ein r. r’ gehe aus r hervor durch eine der folgenden Operationen:
Einfügen von 6 aufeinanderfolgenden "0" oder "1" oder Entfernen von sechs aufeinanderfolgenden "0" oder "1".

Sei $r = \sum_{i=1}^{\infty} z_i 2^{-i}$. Einfügen von sechs "0" ab der n-ten Stelle ergibt $r' = \sum_{i=1}^{n-1} z_i 2^{-i} + \sum_{i=n}^{\infty} z_i 2^{-i-6}$.

Einfügen von sechs "1" ab der n-ten Stelle ergibt $r'' = \sum_{i=1}^{n-1} z_i 2^{-i} + \sum_{i=n}^{n+5} 2^{-i} + \sum_{i=n}^{\infty} z_i 2^{-i-6}$. $r'' = r' + \sum_{i=n}^{n+5} 2^{-i}$.

Definition 5 (Äquivalenz von Wegen)

$r_1$ und $r_2$ heißen äquivalent, wenn $r_1$ durch endlich viele solche Operationen aus $r_2$ erhalten werden kann (in beliebiger Reihenfolge), die Wege $W(r_1)$, $W(r_2)$ heißen ebenfalls äquivalent.

Es handelt sich um eine Äquivalenzrelation.

Die Differenz zweier äquivalenter Zahlen ist rational.

Sei $r_1 = \sum_{i=1}^{n} z_i 2^{-i}$, $z_n = 1$. Sei $r_2 = \sum_{i=1}^{n-1} z_i 2^{-i} + \sum_{i=n+1}^{\infty} 2^{-i}$, so ist $r_1 = r_2$.
$W(r_1)$ bzw. $W(r_2)$ unterscheiden sich nur dadurch, daß das 6-Eck, das ab $z_{n+1}$ durchlaufen wird, linksherum bzw. rechtsherum durchlaufen wird.

Seien r, r' äquivalent und r ≠ r'. Ihre äquivalenten Wege unterschieden sich nur in einem Anfangsteil endlicher Länge, d. h. von einem Punkt A im Gitter an sind sie gleich, die Teilwege von S bis A sind verschieden und i. a. verschieden lang.
Definition der Klasse K(M)
Ein Weg W gehört zur Klasse K(M) wenn für alle Punkte A des Weges W gilt d(A) < M.

Fragen

Für welche Zahlen r bleibt der Weg W(r) in einem beschränkten Gebiet?
Welche Zahlenwege bleiben in einem Winkelraum?
Welche Zahlen haben Wege, die sich nicht selbst kreuzen (kein Gitterpunkt wird zweimal erreicht)?
Für welche Zahlen durchlaufen die Wege alle Gitterpunkte genau n mal, mindestens n mal, n=1, 2, ... ?

Für welche Zahlen r gilt: es gibt zu jedem M Gitterpunkte G, G' ∈ W(r) sodaß (t(W, G) > M, t(W, G') < - M.

Gibt es in jeder Äquivalenzklasse zu jedem M ein r (r(M)) mit W(r) ∉ K(M) ?

Für welche Zahlen r und deren Wege W(r) gilt: zu jedem N gibt es G ∈ W(r) sodaß w(W, G) > N und G' ∈ W(r) sodaß w(W, G') < -N.

Für welche Zahlen r und deren Wege W(r) gilt: es gibt K>0 und zu jedem N gibt es $G_i, G_j$ ∈ W(r) sodaß d($G_i$) > N und $G_j$ mit d($G_j$) < K und i < j.
d. h. der Weg geht beliebig weit weg von S und kommt immer wieder "in die Nähe" von S.

**Verallgemeinerungen**

Wenn man für "0" und "1" links und rechts vertauscht, erhält man einen symetrischen Weg.

Wenn wir unser Gitter mit Quadraten bzw. mit gleichseitigen Dreiecken definieren und r zur Basis 3 bzw. 5 darstellen, können wir analog Wege definieren. Wie hängen diese Wege ab davon wie die Richtungsänderungen zu den Ziffern zugeordnet werden?

In einem quadratischen Gitter sei "0" nach "links", "1" "geradeaus" und "2" nach "rechts". Dann ergeben vier aufeinanderfolgende "0" bzw. "2" eine geschlossenen Teilweg linksherum bzw. rechtsherum.
In Einem Gitter aus gleichseitigen Dreiecken mit analoger Zuordnung zu "0" bis "4" ergeben drei "0" und drei "4" geschlossene Teilwege und sechs "1" und sechs "3" ebenfalls geschlossenen Teilwege.

Welche Eigenschaften der Wege sind von der Basis unabhängig ?

Sind die Eigenschaften der Wege ( Windungszahl, beschränkt, ... ) abhängig von der Basis? Für welche Zahlen ist ein beschränkter Weg W(r) von r zur Basis 2 auch als Weg von r zur Basis 3 (oder 5) beschränkt?

Weitere Verallgemeinerungen, z.B. in höhere Dimensionen liegen nahe.

Um zu einer Zahl r zur Basis 5 einen Weg in einem dreidimensionalen, kubischen Gitter zu definieren, muß zusätzlich zur Anfangsrichtung eine Anfangsorientierung definiert werden, damit zu den Ziffern die Richtungsänderung festgelegt werden kann.

Zu welchen Basen gibt es passende Gitter, in welcher Dimension?

Lorenz Friess
Isnyerstr. 14
D-89079 Ulm
math@lorenz-friess.de